\documentclass[11pt,reqno]{amsart}
\usepackage{amscd,amssymb,amsmath,amsthm}
\usepackage[arrow,matrix]{xy}
\usepackage{graphicx}
\usepackage{cite}
\numberwithin{equation}{section} 

\newcommand{\bc}{\begin{center}}
\newcommand{\ec}{\end{center}}
\newcommand{\be}{\begin{equation}}
\newcommand{\ee}{\end{equation}}
\newcommand{\bea}{\begin{eqnarray}}
\newcommand{\eea}{\end{eqnarray}}
\newcommand{\ba}{\begin{array}}
\newcommand{\ea}{\end{array}}

\begin{document}

\title[$p$-Adic $(3,2)$-rational dynamical systems with three fixed points]{$p$-Adic $(3,2)$-rational dynamical systems with three fixed points}

\author{I.A.Sattarov}

 \address{I.\ A.\ Sattarov \\ Institute of mathematics,
81, Mirzo Ulug'bek str., 100170, Tashkent, Uzbekistan.} \email
{sattarovi-a@yandex.ru}

\begin{abstract}

In this paper we consider dynamical systems generated by $(3,2)$-rational functions
on the field of $p$-adic complex numbers. Each such function has three fixed points.
 We show that Siegel disks of the dynamical system may either coincide or be disjoint
 for different fixed points.
Also, we find the basin of each attractor of the dynamical system.
We show that, for some values of the parameters, there are trajectories
which go arbitrary far from the fixed points.

\end{abstract}

\keywords{Rational dynamical systems; attractors; Siegel disk;
complex $p$-adic field.} \subjclass[2010]{46S10, 12J12, 11S99,
30D05, 54H20.} \maketitle

\section{Introduction}

It is known that in $p$-adic analysis, rational functions play a role
similar to the role of analytic functions in complex analysis.
Therefore, it is natural to study the dynamical systems generated by
these functions in the field of $p$-adic numbers.
Moreover, these $p$-adic dynamical systems appear in the study of
non Archimedean models of physics and biology (see for example \cite{ARS}-\cite{MR2}, \cite{R}).

Let $Q$ be the field of rational numbers.
The completion of $Q$ with  respect to $p$-adic norm defines the
$p$-adic field which is denoted by $Q_p$.

The algebraic completion of $Q_p$ is denoted by $C_p$ and it is
called {\it complex $p$-adic numbers}.  For any $a\in C_p$ and
$r>0$ denote
$$
U_r(a)=\{x\in C_p : |x-a|_p< r\},\ \ V_r(a)=\{x\in C_p :
|x-a|_p\leq r\},
$$
$$
S_r(a)=\{x\in C_p : |x-a|_p= r\}.
$$

A function $f:U_r(a)\to C_p$ is said to be {\it analytic} if it
can be represented by
$$
f(x)=\sum_{n=0}^{\infty}f_n(x-a)^n, \ \ \ f_n\in C_p,
$$ which converges uniformly on the ball $U_r(a)$.

Now let $f:U\to U$ be an analytic function. Denote $x_n=f^n(x_0)$,
where $x_0\in U$ and $f^n(x)=\underbrace{f\circ\dots\circ
f}_n(x)$.

Recall some  the standard terminology of the theory of dynamical
systems (see for example \cite{PJS}). If $f(x_0)=x_0$ then $x_0$
is called a {\it fixed point}. The set of all fixed points of $f$
is denoted by Fix$(f)$. A fixed point $x_0$ is called an {\it
attractor} if there exists a neighborhood $V(x_0)$ of $x_0$ such
that for all points $y\in V(x_0)$ it holds
$\lim\limits_{n\to\infty}y_n=x_0$. If $x_0$ is an attractor then
its {\it basin of attraction} is
$$
A(x_0)=\{y\in C_p :\ y_n\to x_0, \ n\to\infty\}.
$$
A fixed point $x_0$ is called {\it repeller} if there  exists a
neighborhood $V(x_0)$ of $x_0$ such that $|f(x)-x_0|_p>|x-x_0|_p$
for $x\in V(x_0)$, $x\neq x_0$. Let $x_0$ be a fixed point of a
function $f(x)$. The ball $V_r(x_0)$ (contained in $U$) is said to
be a {\it Siegel disk} if each sphere $S_{\rho}(x_0)$, $\rho<r$ is an
invariant sphere of $f(x)$, i.e. if $x\in S_{\rho}(x_0)$ then all
iterated points $x_n\in S_{\rho}(x_0)$ for all $n=1,2\dots$.  The
union of all Siegel desks with the center at $x_0$ is said to {\it
a maximum Siegel disk} and is denoted by $SI(x_0)$.

 Let $x_0$ be a fixed point of an analytic function  $f(x)$. Put
$$
\lambda=\frac{d}{dx}f(x_0).
$$

The point $x_0$ is {\it attractive} if $0\leq |\lambda|_p<1$, {\it
indifferent} if $|\lambda|_p=1$, and {\it repelling} if $|\lambda|_p>1$.

A function is called a $(n,m)$-rational function if and only if it
can be written in the form $f(x)={P_n(x)\over T_m(x)}$, where
$P_n(x)$ and $T_m(x)$ are polynomial functions with degree $n$ and
$m$ respectively, $T_m(x)$ is non zero polynomial (see
\cite{ARS}).

Since the behavior of the dynamical system depends on the set of fixed points,
when the number of fixed points is fixed, each case has its own character of dynamics.
In recent paper \cite{RS} the dynamical systems of the
$p$-adic $(3,2)$-rational functions with unique fixed point were studied.
In paper \cite{S} the case of two fixed points was considered.
In this paper we investigate the behavior of trajectories of the
dynamical system in the case of three fixed points.

\section{$(3,2)$-Rational dynamical systems with three fixed points}

In this paper we consider the dynamical system associated with the
$(3,2)$-rational function $f:C_p\to C_p$ defined by
\begin{equation}\label{fa}
f(x)=ax\left(\frac{x+b}{x+c}\right)^2, \ \ a(a-1)(b-c)(ab^2-c^2)\ne 0, \ \  a,b,c\in C_p.
\end{equation}
where  $x\neq \hat x=-c$.

\textbf{Remark 1.} We note that if $a=0$ or $b=c$, then $f$ isn't $(3,2)$-rational function.
If $a=1$ or $ab^2=c^2$, then function (\ref{fa}) hasn't three distinct fixed points. Therefore we assumed $a(a-1)(b-c)(ab^2-c^2)\ne 0$.

Note that, function (\ref{fa}) has three fixed points $x_0=0$, $$x_1=-\frac{b\sqrt{a}-c}{\sqrt{a}-1} \ \ \mbox{and} \ \ x_2=-\frac{b\sqrt{a}+c}{\sqrt{a}+1}.$$

So we have $f'(x_0)={{ab^2}\over {c^2}}$, $$f'(x_1)=1+{{2(c-b\sqrt{a})(\sqrt{a}-1)}\over{(c-b)\sqrt{a}}} \ \ {\rm and} \ \ f'(x_2)=1+{{2(c+b\sqrt{a})(\sqrt{a}+1)}\over{(c-b)\sqrt{a}}}.$$

For any $x\in C_p$, $x\ne \hat x$, by simple calculations we
get
\begin{equation}\label{f2}
    |f(x)|_p=|a|_p|x|_p\cdot{{|x+b|_p^2}\over {|x+c|_p^2}}.
\end{equation}

Denote
$\mathcal P=\{x\in C_p: \exists n\in N\cup\{0\}, f^n(x)=\hat x$\}.

By using (\ref{f2}) we define the following functions

1. For $|b|_p<|c|_p$ define the function $\varphi_{a,b,c}:
[0,+\infty)\to [0,+\infty)$ by
$$\varphi_{a,b,c}(r)=\left\{\begin{array}{lllll}
{{|ab^2|_p}\over{|c^2|_p}}r, \ \ {\rm if} \ \ r<|b|_p\\[2mm]
b^*, \ \ \ \ \ \ \ \ {\rm if} \ \ r=|b|_p\\[2mm]
{{|a|_p}\over {|c^2|_p}}r^3, \ \ {\rm if} \ \ |b|_p<r<|c|_p\\[2mm]
c^*, \ \ \ \ \ \ \ \ {\rm if} \ \ r=|c|_p\\[2mm]
|a|_pr, \ \ \ \ \ {\rm if} \ \ r>|c|_p
\end{array}
\right.
$$
where $b^*$ and $c^*$ some positive numbers with
$b^*\leq{{|ab^3|_p}\over{|c^2|_p}}$, $c^*\geq |ac|_p$.

2. For $|b|_p=|c|_p$ define the function $\phi_{a,b}:
[0,+\infty)\to [0,+\infty)$ by
$$\phi_{a,b}(r)=\left\{\begin{array}{ll}
|a|_pr, \ \ {\rm if} \ \ r\neq|b|_p\\[2mm]
\hat b, \ \ \ \ \ \ \ {\rm if} \ \ r=|b|_p
\end{array}
\right.
$$
where  $\hat b$ some positive number.

3. For $|b|_p>|c|_p$ define the function $\psi_{a,b,c}:
[0,+\infty)\to [0,+\infty)$ by
$$\psi_{a,b,c}(r)=\left\{\begin{array}{lllll}
{{|ab^2|_p}\over{|c^2|_p}}r, \ \ {\rm if} \ \ r<|c|_p\\[2mm]
c', \ \ \ \ \ \ \ \ {\rm if} \ \ r=|c|_p\\[2mm]
{{|ab^2|_p}\over{r}}, \ \ \ {\rm if} \ \ |c|_p<r<|b|_p\\[2mm]
b', \ \ \ \ \ \ \ \ {\rm if} \ \ r=|b|_p\\[2mm]
|a|_pr, \ \ \ \ {\rm if} \ \ r>|b|_p
\end{array}
\right.
$$
where $b'$ and $c'$ some positive numbers with
$b'\leq |ab|_p$, $c'\geq {{|ab^2|_p}\over{|c|_p}}$.

Using the formula (\ref{f2}) we easily get the following:

\textbf{Lemma 1.} If $x\in S_r(0)$, $x\neq\hat x$ then the following formula
holds for function (\ref{fa})
$$|f^n(x)|_p=\left\{\begin{array}{lll}
\varphi_{a,b,c}^n(r), \ \ \mbox{if} \ \ |b|_p<|c|_p\\[2mm]
\phi_{a,b}^n(r), \ \ \ \ \mbox{if} \ \ |b|_p=|c|_p\\[2mm]
\psi_{a,b,c}^n(r), \ \ \mbox{if} \ \ |b|_p>|c|_p.
\end{array}\right.$$

Thus the $p$-adic dynamical system $f^n(x), n\geq 1, x\in C_p, x\ne \hat x$ is
related to the real dynamical systems generated by $\varphi_{a,b,c}$, $\phi_{a,b}$ and $\psi_{a,b,c}$. Now we are going to study these (real) dynamical systems.

\textbf{Lemma 2.} The dynamical system generated by $\varphi_{a,b,c}(r), |b|_p<|c|_p$ has the following properties:
\begin{itemize}
\item[1.] ${\rm Fix}(\varphi_{a,b,c})=\{0\}\cup$
$$\left\{\begin{array}{lllll}
\{|c|_p:\, {\rm if}\, c^*=|c|_p\}, \ \ \ \ \ \ \ \ \ \ \ \ \ \ \ \ \ \ \ \ \, {\rm for} \, |a|_p<1,\\[2mm]
\{r: r>|c|_p\} \cup\{|c|_p:\, {\rm if}\, c^*=|c|_p\}, \, {\rm for} \, |a|_p=1,\\[2mm]
\{{{|c|_p}\over{\sqrt{|a|_p}}}\}, \ \ \ \ \ \ \ \ \ \ \ \ \ \ \ \ \ \ \ \ \ \ \ \ \ \ \ \ \ \ \ \ \ \ {\rm for} \, |a|_p>1, \, |ab^2|_p<|c^2|_p,\\[2mm]
\{r: r<|b|_p\}\cup\{|b|_p:\, {\rm if}\, b^*=|b|_p\} , \, {\rm for} \,|ab^2|_p=|c^2|_p,\\[2mm]
\{|b|_p:\, {\rm if}\, b^*=|b|_p\}, \ \ \ \ \ \ \ \ \ \ \ \ \ \ \ \ \ \ \ \ \, {\rm for} \,|ab^2|_p>|c^2|_p.
\end{array}
\right.$$

\item[2.] For $|a|_p<1$, we have
\begin{itemize}
 \item[2.1)] If $r\not\in B=\{|a|^{-k}_p|c|_p: k=0,1,2,\dots\}$, then
$$\lim_{n\to\infty}\varphi_{a,b,c}^n(r)=0.$$

\item[2.2)] If $r\in B$ and $c^*\in B$, then there exists $k\geq 0$, such that\\ $c^*=|a|^{-k}_p|c|_p$ and the sequence
$\mathcal C=\{|a|^{-i}_p|c|_p: \, i=0,1,...,k \}$
is a $(k+1)$-cycle of $\varphi_{a,b,c}$.

\item[2.3)] If $r\in B$ and $c^*\notin B$, then
$$\lim_{n\to\infty}\varphi_{a,b,c}^n(r)=0.$$
\end{itemize}

\item[3.] For $|a|_p=1$, we have
$$\lim_{n\to\infty}\varphi_{a,b,c}^n(r)=\left\{\begin{array}{lll}
0, \ \ \mbox{for all} \ \ r<|c|_p,\\[2mm]
r, \ \ \mbox{for all} \ \ r>|c|_p,\\[2mm]
c^*, \ \ \mbox{if} \ \ r=|c|_p
\end{array}\right.;
$$

\item[4.] If $|a|_p>1$, then
\begin{itemize}
\item[4.1)]
If $|ab^2|_p<|c^2|_p$, then $$\lim_{n\to\infty}\varphi_{\alpha,\beta}^n(r)=\left\{\begin{array}{lll}
0, \ \ \mbox{for all} \ \ r<{{|c|_p}\over{\sqrt{|a|_p}}},\\[2mm]
r, \ \ \mbox{for} \ \ r={{|c|_p}\over{\sqrt{|a|_p}}},\\[2mm]
+\infty, \ \ \mbox{if} \ \ r>{{|c|_p}\over{\sqrt{|a|_p}}}
\end{array}\right.;
$$

\item[4.2)]
If $|ab^2|_p=|c^2|_p$, then $$\lim_{n\to\infty}\varphi_{a,b,c}^n(r)=\left\{\begin{array}{lll}
r, \ \ \mbox{for all} \ \ r<|b|_p,\\[2mm]
b^*, \ \ \mbox{for} \ \ r=|b|_p,\\[2mm]
+\infty, \ \ \mbox{if} \ \ r>|b|_p
\end{array}\right.;
$$

 \item[4.3)] If $|ab^2|_p>|c^2|_p$ and $r\not\in L=\{|a^{-k}b^{-2k+1}c^{2k}|_p: k=0,1,2,\dots\}$, then
$$\lim_{n\to\infty}\varphi_{a,b,c}^n(r)=+\infty.$$

\item[4.4)] If $|ab^2|_p>|c^2|_p$, $r\in L$ and $b^*\in L$, then there exists $k\geq 0$, such that $b^*=|a^{-k}b^{-2k+1}c^{2k}|_p$ and the sequence
$\mathcal D=\{|a^{-i}b^{-2i+1}c^{2i}|_p: \, i=0,1,...,k \}$
is a $(k+1)$-cycle of $\varphi_{a,b,c}$.

\item[4.5)] If $|ab^2|_p>|c^2|_p$, $r\in L$ and $b^*\notin L$, then
$$\lim_{n\to\infty}\varphi_{a,b,c}^n(r)=+\infty.$$
\end{itemize}
\end{itemize}

\textbf{Proof.} 1. This is the result of a simple analysis of the equation $\varphi_{a,b,c}(r)=r$.

Proofs of parts 2-4 follow from the property that $\varphi_{a,b,c}(r)$, $r\notin \{|b|_p,|c|_p\}$ is an increasing
function.

\textbf{Lemma 3.} The dynamical system generated by $\phi_{a,b}(r)$ has the following properties:
\begin{itemize}
\item[A.] ${\rm Fix}(\phi_{a,b})=\{0\}\cup$
$$\left\{\begin{array}{ll}
\{|b|_p:\, {\rm if}\, |b|_p=\hat b\}, \, {\rm for} \, |a|_p\ne 1,\\[2mm]
\{r: \, {\rm if}\, r\ne |b|_p\}\cup\{|b|_p:\, {\rm if}\, |b|_p=\hat b\},\, {\rm for} \,|a|_p=1.
\end{array}
\right.;$$

\item[B.] For $|a|_p<1$, we have
\begin{itemize}
 \item[B.a)] If $r\not\in H=\{|a^{-k}b|_p: k=0,1,2,\dots\}$, then
$$\lim_{n\to\infty}\phi_{a,b}^n(r)=0.$$

\item[B.b)] If $r\in H$ and $\hat b\in H$, then there exists $k\geq 0$, such that $\hat b=|a^{-k}b|_p$ and the sequence
$\mathcal F=\{|a^{-i}b|_p: \, i=0,1,...,k \}$
is a $(k+1)$-cycle of $\phi_{a,b}$.

\item[B.c)] If $r\in H$ and $\hat b\notin H$, then
$$\lim_{n\to\infty}\phi_{a,b}^n(r)=0.$$
\end{itemize}
\item[C.] If $|a|_p=1$, then
$$\lim_{n\to\infty}\phi_{a,b}^n(r)=\left\{\begin{array}{ll}
\hat b, \ \ \mbox{if}\ \ r=|b|_p,\\[2mm]
r,\ \ \mbox{for any}\ \ r\neq |b|_p;
\end{array}
\right.$$

\item[D.] For $|a|_p>1$, we have
\begin{itemize}
 \item[D.a)] If $r\not\in H$, then
$$\lim_{n\to\infty}\phi_{a,b}^n(r)=+\infty.$$

\item[D.b)] If $r\in H$ and $\hat b\in H$, then there exists $k\geq 0$, such that\\
 $\hat b=|a^{-k}b|_p$ and the sequence
$\mathcal F$
is a $(k+1)$-cycle of $\phi_{a,b}$.

\item[D.c)] If $r\in H$ and $\hat b\notin H$, then
$$\lim_{n\to\infty}\phi_{a,b}^n(r)=+\infty.$$
\end{itemize}
\end{itemize}

\textbf{Proof.} Since $\phi_{a,b}(r)$ is a piecewise linear function the proof consists of simple computations, using the graph of the function and varying the parameters $a, \, b$.

The following lemma is obvious:

\textbf{Lemma 4.} The dynamical system generated by $\psi_{a,b,c}(r),\, |b|_p>|c|_p$ has the following properties:
\begin{itemize}
\item[(I)] ${\rm Fix}(\psi_{a,b,c})=\{0\}\cup$
$$\left\{\begin{array}{lllll}
\{|c|_p:\, {\rm if}\, c'=|c|_p\}, \ \ \ \ \ \ \ \ \ \ \ \ \ \ \ \ \ \ \ \ {\rm for} \, |a|_p<1, \, |ab^2|_p<|c^2|_p,\\[2mm]
\{r:\, r<|c|_p\}\cup\{|c|_p:\, {\rm if}\, c'=|c|_p\}, \, {\rm for} \, |a|_p<1, \, |ab^2|_p=|c^2|_p,\\[2mm]
\{|b|_p\sqrt{|a|_p}\}, \ \ \ \ \ \ \ \ \ \ \ \ \ \ \ \ \ \ \ \ \ \ \ \ \ \ \ {\rm for} \, |a|_p<1, \, |ab^2|_p>|c^2|_p,\\[2mm]
\{r:\, r>|b|_p\}\cup\{|b|_p:\, {\rm if}\, b'=|b|_p\}, \, {\rm for} \, |a|_p=1,\\[2mm]
\{|b|_p: \, {\rm if}\, b'=|b|_p\},\ \ \ \ \ \ \ \ \ \ \ \ \ \ \ \ \ \ \ \ \ {\rm for} \,|a|_p>1.
\end{array}
\right.;$$

\item[(II)] Let $|a|_p<1$ and $|ab^2|_p<|c^2|_p$.
\begin{itemize}
 \item[II.i)] If $r\not\in B=\{|a|^{-k}_p|c|_p: k=0,1,2,\dots\}$, then
$$\lim_{n\to\infty}\psi_{a,b,c}^n(r)=0.$$

\item[II.ii)] If $r\in B$ and $c'\in B$, then there exists $k\geq 0$, such that\\ $c'=|a|^{-k}_p|c|_p$ and the sequence
$\mathcal C=\{|a|^{-i}_p|c|_p: \, i=0,1,...,k \}$
is a $(k+1)$-cycle of $\psi_{a,b,c}$.

\item[II.iii)] If $r\in B$ and $c'\notin B$, then
$$\lim_{n\to\infty}\psi_{a,b,c}^n(r)=0.$$
\end{itemize}
\item[(III)]  Let $|a|_p<1$ and $|ab^2|_p=|c^2|_p$.
\begin{itemize}
 \item[III.i)] If $r\not\in B$, then there exists an integer $k\geq 0$, such that \\
$\psi_{a,b,c}^n(r)=\psi_{a,b,c}^k(r),\ \ \mbox{for any} \ \ n\geq k.$

\item[III.ii)] If $r\in B$ and $c'\in B$, then there exists $k\geq 0$, such that\\ $c'=|a|^{-k}_p|c|_p$ and the sequence
$\mathcal C=\{|a|^{-i}_p|c|_p: \, i=0,1,...,k \}$
is a $(k+1)$-cycle of $\psi_{a,b,c}$.

\item[III.iii)] If $r\in B$ and $c'\notin B$, then
there exists an integer $k\geq 0$, such that
$\psi_{a,b,c}^n(r)=\psi_{a,b,c}^k(r),\ \ \mbox{for any} \ \ n\geq k.$
\end{itemize}
\item[(IV)] If $|a|_p<1$ and $|ab^2|_p>|c^2|_p$, then there exists the invariant set $\Lambda=(|b|_p\sqrt{|a|_p}-\lambda, |b|_p\sqrt{|a|_p}+\lambda)$, such that $\psi_{a,b,c}^2(r)=r$ for any $r\in\Lambda$, moreover, if $r\notin\Lambda\cup\{0\}$, then there exists $k\geq 1$, such that $\psi_{a,b,c}^k(r)\in\Lambda$,\\
     where $\lambda=\min\{|b|_p\sqrt{|a|_p}-|c|_p, |b|_p(1-\sqrt{|a|_p})\}$.
\item[(V)] Let $|a|_p=1$.
\begin{itemize}
 \item[V.i)] If $r\not\in L=\{|a^{-k}b^{-2k+1}c^{2k}|_p: k=0,1,2,\dots\}$, then there exists $k\geq 0$, such that
$\psi_{a,b,c}^n(r)=\psi_{a,b,c}^k(r)$, for any $n\geq k.$

\item[V.ii)] If $r\in L$ and $b'\in L$, then there exists $k\geq 0$, such that\\ $b'=|a^{-k}b^{-2k+1}c^{2k}|_p$ and the sequence
$\mathcal D$ is a $(k+1)$-cycle of $\psi_{a,b,c}$.
\item[V.iii)] If $r\in L$ and $b'\notin L$, then
there exists an integer $k\geq 0$, such that
$\psi_{a,b,c}^n(r)=\psi_{a,b,c}^k(r),\ \ \mbox{for any} \ \ n\geq k.$
\end{itemize}

\item[VI.] For $|a|_p>1$, we have
\begin{itemize}
 \item[VI.i)] If $r\not\in L$, then
$$\lim_{n\to\infty}\psi_{a,b,c}^n(r)=+\infty.$$

\item[VI.ii)] If $r\in L$ and $b'\in L$, then there exists $k\geq 0$, such that\\
$b'=|a^{-k}b^{-2k+1}c^{2k}|_p$ and the sequence
$\mathcal D$
is a $(k+1)$-cycle of $\psi_{a,b,c}$.

\item[VI.iii)] If $r\in L$ and $b'\notin L$, then
$$\lim_{n\to\infty}\psi_{a,b,c}^n(r)=+\infty.$$
\end{itemize}
\end{itemize}

Now we shall apply these lemmas to the study of the $p$-adic dynamical system generated by $f$.

For $x\in S_{|b|_p}(0)$, we denote
$$b^*(x)=|a|_p|b|_p\cdot{{|x+b|_p^2}\over {|x+c|_p^2}}.$$

For $x\in S_{|c|_p}(0)$, we denote
$$c^*(x)=|a|_p|c|_p\cdot{{|x+b|_p^2}\over {|x+c|_p^2}}.$$
Using Lemma 1 and Lemma 2 we obtain the following

\textbf{Theorem 1.} If $|b|_p<|c|_p$ and $x\in S_r(0)$, then
 the $p$-adic dynamical system generated by $f$ has the following properties:
\begin{itemize}
\item[1.] The following spheres are invariant with respect to $f$:
$$\begin{array}{lll}
S_r(0),\, {\rm if} \,r>|c|_p,\,  |a|_p=1,\\[2mm]
S_r(0), \, {\rm if} \, r<|b|_p, \, |ab^2|_p=|c^2|_p,\\[2mm]
S_{{|c|_p}\over{\sqrt{|a|_p}}}(0), \, {\rm if}\, |ab^2|_p<|c^2|_p, \, |a|_p>1.\\[2mm]
\end{array}
;$$

\item[2.] For $|a|_p<1$, we have
\begin{itemize}
\item[2.1)] If $r\notin B=\{|a|^{-k}_p|c|_p: k=0,1,2,\dots\}$, then $$\lim_{n\to\infty}f^n(x)=0.$$
\item[2.2)] If $r\in B$, then there exists $k\geq 0$, such that $f^k(x)\in S_{|c|_p}(0)$ and if $c^*(f^k(x))\notin B$, then $$\lim_{n\to\infty}f^n(x)=0.$$
\item[2.3)] $x_1, \, x_2\in S_{|c|_p}(0)$ and $$|x_1-x_2|_p=\left\{\begin{array}{ll}
{{|c|_2}\over 2}, \ \ \mbox{for} \ \ p=2,\\[2mm]
|c|_p, \ \ \mbox{for} \ \ p\geq3.
\end{array}\right.$$
\item[2.4)] If $p\geq 3$, then the fixed points $x_1$ and $x_2$ are repeller and the inequality $|f(x)-x_i|_p>|x-x_i|_p$ is hold for all $x\in U_{|c|_p}(x_i)$, $i=1,2.$
\item[2.5)] Let $p=2$. Then we have the following:

if $|a|_2<{1\over 4}$, then the fixed points $x_1$ and $x_2$ are repeller.

if $|a|_2={1\over 4}$, then the fixed points $x_1$ and $x_2$ are attractor.

if $|a|_2>{1\over 4}$, then the fixed points $x_1$ and $x_2$ are indifferent.
\end{itemize}
\item[3.] If $|a|_p=1$, then $A(0)=U_{|c|_p}(0)$. Moreover,
$x_i\notin V_{|c|_p}(0)$, the fixed point $x_i$ is indifferent if $p=2$ and
the fixed point $x_i$ is indifferent or attractor for $p\geq 3$, $i=1, \, 2$.
\item[4.] Let $|a|_p>1$.
\begin{itemize}
\item[4.1)]
If $|ab^2|_p<|c^2|_p$, then $A(0)=U_{{{|c|_p}\over{\sqrt{|a|_p}}}}(0)$ and for all $r>{{|c|_p}\over{\sqrt{|a|_p}}}$ we have $$\lim_{n\to\infty}|f^n(x)|_p=+\infty.$$
Moreover,
$x_i\in S_{{{|c|_p}\over{\sqrt{|a|_p}}}}(0)$, the fixed point $x_i$ is indifferent if $p=2$ and
the fixed point $x_i$ is indifferent or attractor for $p\geq 3$, $i=1, \, 2$.
\item[4.2)]
If $|ab^2|_p=|c^2|_p$, then $SI(0)=U_{|b|_p}(0)$ and for all $r>|b|_p$ we have $$\lim_{n\to\infty}|f^n(x)|_p=+\infty.$$
Moreover,
$x_i\in V_{|b|_p}(0)$, the fixed point $x_i$ is indifferent if $p=2$ and
the fixed point $x_i$ is indifferent or attractor for $p\geq 3$, $i=1, \, 2$.
\item[4.3)]
If $|ab^2|_p>|c^2|_p$ and $r\notin L=\{|a^{-k}b^{-2k+1}c^{2k}|_p: k=0,1,2,\dots\}$, then $$\lim_{n\to\infty}|f^n(x)|_p=+\infty.$$

\item[4.4)]
If $|ab^2|_p>|c^2|_p$ and $r\in L$, then there exists $k\geq 0$, such that $f^k(x)\in S_{|b|_p}(0)$ and if $b^*(f^k(x))\notin L$, then $$\lim_{n\to\infty}|f^n(x)|_p=+\infty.$$

\item[4.5)]
If $|ab^2|_p>|c^2|_p$, then $x_i\in S_{|b|_p}(0)$ and the fixed point $x_i, \, i=1, \, 2$ is repeller for $p\geq 3$.

Moreover, if $p=2$ and $|b|_2\sqrt{|a|_2}=2|c|_2$, then  the fixed point $x_i, \, i=1, \, 2$ is attractor and if $|b|_2\sqrt{|a|_2}>2|c|_2$, then the fixed point $x_i, \, i=1, \, 2$ is repeller.
\end{itemize}
\end{itemize}

By Lemma 1 and Lemma 3 we obtain the following

\textbf{Theorem 2.} If $|b|_p=|c|_p$ and $x\in S_r(0)$, then the $p$-adic dynamical system generated by $f$ has the following properties:
\begin{itemize}

\item[A.] Let $|a|_p<1$. Then:
\begin{itemize}
\item[A.a)] If $r\notin H=\{|a^{-k}b|_p: k=0,1,2,\dots\}$, then $$\lim_{n\to\infty}f^n(x)=0.$$
\item[A.b)] If $r\in H$, then there exists $k\geq 0$ such that $r=|a^{-k}b|_p$ and $f^k(x)\in S_{|b|_p}(0)$.
\item[A.c)] If $x\in S_{|b|_p}(0)$ and $b^*(x)\notin H$ then $$\lim_{n\to\infty}f^n(x)=0.$$
\item[A.d)] If $x\in S_{|b|_p}(0)$ and $b^*(x)\in H$ then there exists $k\geq 0$
such that $b^*(x)=|a^{-k}b|_p$ and $f^{k+1}(x)\in S_{|b|_p}(0)$.
\item[A.e)] The fixed points $x_1$, $x_2$ are repeller and $x_i\in S_{|b|_p}(0), \, i=1,2.$
\end{itemize}
\item[B.] Let $|a|_p=1$. Then the sphere $S_r(0)$ is invariant for $f$ if $r\neq |b|_p$.
If $x\in S_{|b|_p}(0)$, then one of the following two possibilities holds:
\begin{itemize}

\item[B.a)] There exists $k\in N$ and $\mu_k\neq |b|_p$ such that $f^n(x)\in S_{\mu_k}(0)$ for any $n\geq k$ and
$f^m(x)\in S_{|b|_p}(0)$ for any $m\leq k-1$.

\item[B.b)] The trajectory $\{f^k(x), k\geq 1\}$ is a subset of
$S_{|b|_p}(0)$.
\end{itemize}

\item[C.] Let $|a|_p>1$. Then:
\begin{itemize}

\item[C.a)] If $r\notin H$, then $$\lim_{n\to\infty}|f^n(x)|_p=+\infty.$$
\item[C.b)] If $r\in H$, then there exists $k\geq 0$ such that $r=|a^{-k}b|_p$ and $f^k(x)\in S_{|b|_p}(0)$.
\item[C.c)] If $x\in S_{|b|_p}(0)$ and $b^*(x)\notin H$ then $$\lim_{n\to\infty}|f^n(x)|_p=+\infty.$$
\item[C.d)] If $x\in S_{|b|_p}(0)$ and $b^*(x)\in H$ then there exists $k\geq 0$
such that $b^*(x)=|a^{-k}b|_p$ and $f^{k+1}(x)\in S_{|b|_p}(0)$.
\item[C.e)] The fixed points $x_1$, $x_2$ are repeller and $x_i\in S_{|b|_p}(0), \, i=1,2.$
\end{itemize}
\end{itemize}
\textbf{Proof.}
A. The proof of parts A.a) -A.d) comes from the parts A, B of Lemma 3.

A.e) Note that $$x_1=-\frac{b\sqrt{a}-c}{\sqrt{a}-1} \ \ \mbox{and} \ \ x_2=-\frac{b\sqrt{a}+c}{\sqrt{a}+1}.$$
Moreover, $$f'(x_1)=1+{{2(c-b\sqrt{a})(\sqrt{a}-1)}\over{(c-b)\sqrt{a}}} \ \ {\rm and} \ \ f'(x_2)=1+{{2(c+b\sqrt{a})(\sqrt{a}+1)}\over{(c-b)\sqrt{a}}}.$$
Therefore, if $|a|_p<1$ and $|b|_p=|c|_p$, then $|x_1|_p=|x_2|_p=|b|_p$. Moreover, $|f'(x_i)|_p>1$, $i=1, \, 2$.

The proof of part B of this Theorem comes from part C of Lemma 3.

Proof of part C is similar to proof of part A.

By Lemma 1 and Lemma 4 we get

\textbf{Theorem 3.} If $|b|_p>|c|_p$ and $x\in S_r(0)$, then the dynamical system generated by $f$ has the following properties:
\begin{itemize}
\item[I.] The following spheres are invariant:

$S_r(0), \ \ r<|c|_p \ \ \ \ \mbox{for} \ \ |a|_p<1, \ \ |ab^2|_p=|c^2|_p;$

$S_{|b|_p\sqrt{|a|_p}}(0), \ \ \ \ \ \ \, \mbox{for} \ \ |a|_p<1, \ \ |ab^2|_p>|c^2|_p;$

$S_r(0), \ \ r>|b|_p \ \ \ \ \mbox{for} \ \ |a|_p=1.$

\item[II.] Let $|a|_p<1$ and $|ab^2|_p<|c^2|_p$. Then:
\begin{itemize}
\item[II.a)] If $r\notin B$, then $$\lim_{n\to\infty}f^n(x)=0.$$
\item[II.b)] If $r\in B$, then there exists $k\geq 0$ such that $r=|a^{-k}c|_p$ and $f^k(x)\in S_{|c|_p}(0)$.
\item[II.c)] If $x\in S_{|c|_p}(0)$ and $c^*(x)\notin B$ then $$\lim_{n\to\infty}f^n(x)=0.$$
\item[II.d)] If $x\in S_{|c|_p}(0)$ and $c^*(x)\in B$ then there exists $k\geq 0$
such that $c^*(x)=|a^{-k}c|_p$ and $f^{k+1}(x)\in S_{|c|_p}(0)$.
\item[II.e)] The fixed points $x_1$, $x_2$ are repeller for $p\geq 3$ and $x_i\in S_{|c|_p}(0), \, i=1,2.$
Moreover, for $p=2$ the fixed points $x_1$ and $x_2$ are repeller if $|c|_2>2|b|_2\sqrt{|a|_2}$, the fixed points $x_1$ and $x_2$ are attractor if $|c|_2=2|b|_2\sqrt{|a|_2}$.
\end{itemize}
\item[III.] Let $|a|_p<1$ and $|ab^2|_p=|c^2|_p$. Then:
\begin{itemize}
\item[III.a)] If $r\notin B$, then there exists $k\geq 0$ such that $f^n(x)\in S_{\psi_{a,b,c}^k(r)}(0)$ for any $n\geq k$.
\item[III.b)] If $r\in B$, then there exists $k\geq 0$ such that $r=|a^{-k}c|_p$ and $f^k(x)\in S_{|c|_p}(0)$.
\item[III.c)] If $x\in S_{|c|_p}(0)$ and $c^*(x)\notin B$ then there exists $k\geq 1$ such that $f^n(x)\in S_{\psi_{a,b,c}^k(c^*(x))}(0)$ for any $n\geq k$.
\item[III.d)] If $x\in S_{|c|_p}(0)$ and $c^*(x)\in B$ then there exists $k\geq 0$
such that $c^*(x)=|a^{-k}c|_p$ and $f^{k+1}(x)\in S_{|c|_p}(0)$.
\item[III.e)] $x_i\in V_{|c|_p}(0), \, i=1,2.$ Fixed points $x_1$ and $x_2$ of the function $f$ may be an attractor or indifferent fixed point.
\end{itemize}
\item[IV.] If $|a|_p<1$ and $|ab^2|_p>|c^2|_p$, then $f^2(S_r(0))\subset S_r(0)$ for any\\ $r\in \Lambda=(|b|_p\sqrt{|a|_p}-\lambda, |b|_p\sqrt{|a|_p}+\lambda)$,\\ where $\lambda=\min\{|b|_p\sqrt{|a|_p}-|c|_p, |b|_p(1-\sqrt{|a|_p})\}$.

Moreover, if $r\notin\Lambda\cup\{0\}$, then there exists $k\geq 1$ such that $\psi_{a,b,c}^k(r)\in\Lambda$ and $f^k(x)\in S_{\psi_{a,b,c}^k(r)}(0)$.

\item[V.] Let $|a|_p=1$. Then:
\begin{itemize}
\item[V.a)] If $r\notin L$, then there exists $k\geq 0$ such that $f^n(x)\in S_{\psi_{a,b,c}^k(r)}(0)$ for any $n\geq k$.
\item[V.b)] If $r\in L$, then there exists $k\geq 0$ such that $r=|a^{-k}b^{-2k+1}c^{2k}|_p$ and $f^k(x)\in S_{|b|_p}(0)$.
\item[V.c)] If $x\in S_{|b|_p}(0)$ and $b^*(x)\notin L$ then there exists $k\geq 1$ such that $f^n(x)\in S_{\psi_{a,b,c}^k(b^*(x))}(0)$ for any $n\geq k$.
\item[V.d)] If $x\in S_{|b|_p}(0)$ and $b^*(x)\in L$ then there exists $k\geq 0$
such that $b^*(x)=|a^{-k}b^{-2k+1}c^{2k}|_p$ and $f^{k+1}(x)\in S_{|b|_p}(0)$.
\item[V.e)] $x_i\notin V_{|b|_p}(0), \, i=1,2.$ Fixed points $x_1$ and $x_2$ of the function $f$ may be an attractor or indifferent fixed point.
\end{itemize}
\item[VI.] Let $|a|_p>1$. Then:
\begin{itemize}
\item[VI.a)] If $r\notin L$, then $$\lim_{n\to\infty}|f^n(x)|_p=+\infty.$$
\item[VI.b)] If $r\in L$, then there exists $k\geq 0$ such that $r=|a^{-k}b^{-2k+1}c^{2k}|_p$ and $f^k(x)\in S_{|b|_p}(0)$.
\item[VI.c)] If $x\in S_{|b|_p}(0)$ and $b^*(x)\notin L$ then $$\lim_{n\to\infty}|f^n(x)|_p=+\infty.$$
\item[VI.d)] If $x\in S_{|b|_p}(0)$ and $b^*(x)\in L$ then there exists $k\geq 0$
such that $b^*(x)=|a^{-k}b^{-2k+1}c^{2k}|_p$ and $f^k(x)\in S_{|b|_p}(0)$.
\item[VI.e)] The fixed points $x_1$, $x_2$ are repeller and $x_i\in S_{|b|_p}(0), \, i=1,2.$
\end{itemize}
\end{itemize}

\section{Acknowledgments}

The author expresses his deep gratitude to Professor U. A. Rozikov for setting up the problem and for the useful suggestions.

 \end{document}